\documentclass[brazil,english]{article}
\usepackage[T1]{fontenc}
\usepackage[latin9]{inputenc}
\usepackage{geometry}
\usepackage{amsmath}

\makeatletter
\@ifundefined{date}{}{\date{}}
\makeatother

\usepackage{babel}
\begin{document}

\title{Solving Polynomial Equations from Complex Numbers}

\selectlanguage{brazil}%

\author{Ricardo S. Vieira%
\thanks{E-mail: rsvieira@df.ufscar.br.%
} \vspace{0.25cm} \\  \textit{Departamento de Física, Universidade Federal de São Carlos, São Paulo, Brasil}}
\maketitle
\selectlanguage{english}%
\begin{abstract}
We show that a polynomial equation of degree less than 5 and with
real parameters can be solved by regarding the variable in which the
polynomial depends as a complex variable. For do it so, we only have
to separate the real and imaginary parts of the resultant polynomial
and solve them separately. 
\end{abstract}

\section{Solving the Quadratic Equation}

Let we begin by considering the quadratic equation 
\begin{equation}
z^{2}+az+b=0,\label{quadratic}
\end{equation}
where $a$ and $b$ are two real numbers. 

Regarding $z$ as a complex variable, we can put 
\begin{equation}
z=x+iy,\label{z1}
\end{equation}
where $x$ and $y$ are two real numbers and $i=\sqrt{-1}$. By replacing
$z$ from (\ref{z1}) in (\ref{quadratic}) we get, after we had separated
the real part from the imaginary one, the following equations, 
\begin{equation}
\begin{cases}
x^{2}-y^{2}+ax+b=0,\\
2xy+ay=0.
\end{cases}\label{S1}
\end{equation}
 Since these equations must be valid for arbitrary $x$ and $y$ we
get, solving the second equation for $x$, 
\begin{equation}
x=-a/2.
\end{equation}
The substitution of this result onto the first equation (\ref{S1})
enable us to find $y$, which is given by 
\begin{equation}
y=\pm\sqrt{b-a^{2}/4},
\end{equation}
 and, therefore, we find that the solutions of (\ref{quadratic})
are 
\begin{equation}
z=-\frac{a}{2}\pm\sqrt{\frac{a^{2}}{4}-b}.
\end{equation}

\section{Solving the Cubic Equation}

Let we consider now a cubic equation as 
\begin{equation}
w^{3}+\alpha w^{2}+\beta w+\gamma=0,\label{cubic-0}
\end{equation}
 where $\alpha$, $\beta$ and $\gamma$ are real parameters. As it
is well know, this equation can be put into the reduced form 

\begin{equation}
z^{3}+az+b=0,\label{cubic}
\end{equation}
through the change of variable $w=z-\alpha/3$. The parameters $a$
and $b$ are two real numbers as well, since they are given by the
expressions 
\begin{equation}
a=\beta-\alpha^{2}/3,\qquad b=2\alpha^{3}/27-\alpha\beta/3+\gamma.\label{ab}
\end{equation}

Therefore, let we concentrate ourselves with the equation (\ref{cubic}).
If we put $z=x+iy$ in that equation, as we have done in the precedent
section, we will obtain the system of equations 
\begin{equation}
\begin{cases}
x^{3}-3xy^{2}+ax+b=0,\\
y^{3}-3x^{2}y-ay=0,
\end{cases}\label{S2}
\end{equation}
by separating the real and imaginary parts of the resultant equation.
Now, solving the second equation for $x$ and replacing the result
on the first, we get 
\begin{equation}
8x^{3}+2ax-b=0,
\end{equation}
which is a third degree equation as (\ref{cubic}). Therefore, it
appears that a cubic equation cannot be solved by the same technique
utilized to solve the quadratic equation. This, however, it is only
apparent, as we will show in the following.

The reason for this apparent failure is that the combination $z=x+iy$
is not the good one for cubic equations, even though for quadratic
equations it is. In fact, notice that in the case of a quadratic equation
we can write $z$ as $x+y\sqrt{-1}$ and hence the coefficient which
multiplies $y$ is, in this case, one of the two \emph{square} \emph{roots}
of $-1$. Therefore, in the case of a cubic equation we can naturally
try the ansatz%
\footnote{Roots of unity plays a central role in the theory of polynomials.
Indeed, with the help of the roots of unity, and by analyzing the
symmetries of the roots of a polynomial equation, Lagrange developed
a profound theory, by which he explained in details why the methods
proposed up to date for solving polynomial equations worked well for
the quadratic, cubic and quartic equations but failed for equations
of higher degree. As it is know, the theory of Lagrange influenced
Ruffini, Abel and Galois, whose work culminated in the proof of impossibility
for solving a general polynomial equation (of degree greater than
4) by radicals. For more details see \cite{Dehn,Tignol,Pesic}.%
} 
\begin{equation}
z=x+\omega y,
\end{equation}
where $\omega$ is one of the complex the \emph{cubic roots} of $-1$
(except $-1$ itself, of course). Actually, we can verify that $\omega$
could be chosen as one of the cubic roots of unity (except $+1$,
of course) as well. 

Then, let we choose, for instance, 
\begin{equation}
\omega=\frac{1+\sqrt{-3}}{2}.\label{omega}
\end{equation}
 Substituting this value for $\omega$ into the equation (\ref{cubic})
and separating the real and imaginary parts, we get the equations,
\begin{equation}
\begin{cases}
x^{3}-y^{3}+3x^{2}y/2-3xy^{2}/2+ay/2+ax+b=0,\\
3xy^{2}+3x^{2}y+ay=0.
\end{cases}\label{S3}
\end{equation}
and solving the second equation of (\ref{S3}) for $y$,
\begin{equation}
y=-x-\frac{a}{3x}.
\end{equation}
Replacing this value of $y$ into the first equation, we will be led
to the equation 
\begin{equation}
x^{6}-bx^{3}-\frac{a^{3}}{27}=0,
\end{equation}
 which is a quadratic equation on $x^{3}$ and, therefore, has the
solutions 
\begin{equation}
x^{3}=\frac{b}{2}\pm\sqrt{\frac{b^{2}}{4}+\frac{a^{3}}{27}}.
\end{equation}
Consequently, the solutions of a cubic equation as (\ref{cubic})
can be written, after some simplification, as 
\begin{equation}
z=\left(1-\omega\right)\sqrt[3]{\frac{b}{2}\pm\sqrt{\frac{b^{2}}{4}+\frac{a^{3}}{27}}}-\cfrac{\omega a/3}{\sqrt[3]{\frac{b}{2}\pm\sqrt{\frac{b^{2}}{4}+\frac{a^{3}}{27}}}},
\end{equation}
with the roots on this equation chosen appropriately, of course. The
solutions of the general cubic equation (\ref{cubic-0}) follows from
the substitution of the parameters $a$ and $b$ by the expressions
(\ref{ab}) and from the relation $w=z-\alpha/3$.

\section{Solving the Quartic Equation}

Let now we consider the equation 
\begin{equation}
w^{4}+\alpha w^{3}+\beta w^{2}+\gamma w+\delta=0,\label{quartic-0}
\end{equation}
where $\alpha$, $\beta$, $\gamma$ and $\delta$ are four real numbers.
As before, this equation can be put into the reduced form 
\begin{equation}
z^{4}+az^{2}+bz+c=0\label{quartic}
\end{equation}
by the change of variable $w=z-\alpha/4$, and the new parameters
$a$, $b$ and $c$, which are real parameters as well, are given
by the expressions 
\begin{equation}
a=\beta-3\alpha^{2}/8,\qquad b=\alpha^{3}/8-\alpha\beta/2+\gamma,\qquad c=-3\alpha^{4}/256+\alpha^{2}\beta/16-\alpha\gamma/4+\delta.\label{abc}
\end{equation}
From now we will consider only the equation (\ref{quartic}). 

Since the quartic roots of the unity are $+1$, $-1$, $+i$ and $-i$
we can choose here the simplest combination for $z$, that is, 
\begin{equation}
z=x+iy,
\end{equation}
where $x$ and $y$ are real. Then, proceeding as before, we replace
this value of $z$ onto the equation (\ref{quartic}) and we separate
their real and imaginary parts. We will get the system 
\begin{equation}
\begin{cases}
x^{4}+y^{4}-6x^{2}y^{2}+ax^{2}-ay^{2}+bx+c=0,\\
4x^{3}y-4xy^{3}+2axy+by=0.
\end{cases}\label{S4}
\end{equation}
The non-trivial solution of the second equation above is 
\begin{equation}
y^{2}=x^{2}+\frac{b}{4x}+\frac{a}{2},
\end{equation}
and the substitution of this value onto the first equation in (\ref{S4})
provides the equation 
\begin{equation}
x^{6}+\frac{a}{2}x^{4}+\left(\frac{a^{2}}{16}-\frac{c}{4}\right)x^{2}-\frac{b^{2}}{64}=0,
\end{equation}
which is a cubic equation on the variable $x^{2}$ and, therefore,
can be solved as well. Once this equation is solved, the solutions
of the general quartic equation (\ref{quartic-0}) can be obtained
by replacing $a,$ $b$ and $c$ by their expressions (\ref{abc})
and from the relation $w=z-\alpha/4$. Although explicit formulæ can
be written, we will not do so because the expressions are too much
extensive.

Finally, we should point out that for a fifth degree equation this
method does not give any appreciable help, since by separating the
real and imaginary parts of that polynomial we shall get other more
complicated equations. These difficulties, however, were expected
already, since the roots of a general equation of degree 5 (or higher)
can not be expressed by radicals, as was proven long time ago by Ruffini,
Abel and Galois.

\end{document}